\documentclass[11pt]{article}
\usepackage{mathrsfs}
\usepackage[centertags]{amsmath}
\usepackage{amsfonts}
\usepackage{amssymb}
\usepackage{amsthm}
\usepackage{amscd,amssymb,amsthm}
\newcommand{\be}{\begin{equation}}
\newcommand{\ee}{\end{equation}}
\newcommand{\bea}{\begin{eqnarray}}
\newcommand{\eea}{\end{eqnarray}}
\newcommand{\bna}{\begin{eqnarray*}}
\newcommand{\ena}{\end{eqnarray*}}

\newtheorem{thm}{Theorem}[section]

\newtheorem{lemma}[thm]{Lemma}

\theoremstyle{definition}

\theoremstyle{remark}

\numberwithin{equation}{section}
\begin{document}
\title{ \bf Br\"{u}ck Conjecture with hyper-order less than one}
\author{\sc  Guowei Zhang
%\thanks{ \indent
%%\Newline\Indent\Quad The First Author
%%Was Supported By The Nnsf Of China (No.10371065), The Nsf Of
%%Shangdong Province, China (No Z2002A01) And The Nsfc-Rfbr
%}
 \\
School of Mathematics and Statistics, Anyang Normal University,
Anyang, Henan\\ 455000,  China, E-mail: herrzgw@foxmail.com }
\date{}
\maketitle

\begin{center}{\bf Abstract}\end{center}

In this paper we affirm Br\"{u}ck conjecture provided $f$ is of hyper-order less than one by studying the infinite hyper-order of solutions of a complex differential equation.\\


\noindent{\bf 2000 MR subject classification}: 30D35; 34M10.

\baselineskip14pt \setcounter{section}{0}

\section{Introduction and main results}

In this article, we assume the reader is familiar with standard
notations and basic results of Nevanlinna's value distribution
theory in the complex plane $\mathbb{C}$, see
\cite{Laine,YY}. The order and hyper-order of an entire function $f$ are defined as

\bea &&\rho(f)=\limsup_{r\rightarrow+\infty}\frac{\log^+ T(r,f)}{\log r}=\limsup_{r\rightarrow+\infty}\frac{\log^+ \log^+M(r,f)}{\log r},\nonumber\\
&&\rho_2(f)=\limsup_{r\rightarrow+\infty}\frac{\log^+\log^+ T(r,f)}{\log r}=\limsup_{r\rightarrow+\infty}\frac{\log^+\log^+ \log^+M(r,f)}{\log r},\nonumber
\eea
respectively, where $M(r,f)$ denotes the maximum modulus of $f$ on the circle $|z|=r$.

If $f, g$ are two meromorphic funcitons in the complex plane, we say $f, g$ share a constant $a$ CM if $f-a$ and $g-a$ have the same zeros with the same multiplicities. Rubel and Yang \cite{RubelYang} proved for a nonconstant entire function, if $f$ and its derivative $f'$ share two finite distinct values CM, then $f\equiv f'$. Later on, Br\"{u}ck \cite{Bruck} constructed entire functions with integer or infinite hyper-order to show that $f$ and $f'$ share 1 CM fails to obtain $f\equiv f'$. Therefore, Br\"{u}ck proposed the following conjecture.\\

\noindent {\bf Br\"{u}ck Conjecture}\cite{Bruck}
{\it Let $f$ be a nonconstant entire function such that its hyper-order is finite but not a positive integer. If $f$ and $f'$ share one finite value $a$ CM, then $f'-a=c(f-a)$, where $c$ is a nonzero constant.}

The conjecture is false in general for meromorphic function $f$, see a counterpart in \cite{GundersenYang}. Br\"{u}ck \cite{Bruck} showed the conjecture is right for the case $a=0$. Later, Gundersen and Yang \cite{GundersenYang} proved the conjecture is true for the case that $f$ is of finite order. Further on, Chen and Shon \cite{ChenShonTW} showed that the conjecture is also true under the condition $f$ is of hyper-order strictly less than 1/2. Recently, Cao \cite{CaoAus} gave an affirmative answer to this conjecture under hypothesis that the hyper-order of $f$ is equal to 1/2.

There are many results closely related to Br\"{u}ck conjecture, mainly in two directions. One is generalizing the shared value $a$ to a nonconstant function, such as polynomial, entire small function respect to $f$, or entire functions with lower order than $f$, (e. g. see \cite{CaoPolon,Changjm,LiXiaomin,LiXMGao,WangJun}) and another is improving the first derivative of $f$ to arbitrary $k$-th derivative (e. g. see \cite{CaoPolon,ChenShonKR,ChenZhang,LiXMGao,YangLZKodai}).
The main purpose of this paper is to confirm Br\"{u}ck conjecture provided the hyper-order of $f$ is less than one. In fact, we obtain the following result.
\begin{thm}\label{T1}  {\it Let $f$ be a nonconstant entire function with hyper-order $\rho_2(f)<1$. If $f$ shares
one finite value $a$ CM with its $k$-th derivative, then $f^{(k)}-a=c(f-a)$, where $c$ is a nonzero constant.}
\end{thm}
 Thus, for the final solving of this conjecture one should consider the remaining case that the hyper-order of $f$ lying in $(1,+\infty)\backslash \mathbb{N}$.
 In order to study this conjecture, many authors paid attention to the nonhomogeneous linear complex differential equation
 \bea f^{(k)}-e^{p(z)}f=Q(z), \ k\in \mathbb{N},\eea
where $p(z)$ is an entire function and $Q(z)$ is a constant or an entire function. Under suitable conditions on $p(z)$ and $Q(z)$ thought of ways to prove the nontrivial solution of this equation is of infinite order or infinite hyper-order, such as \cite{CaoPolon,ChenShonTW,ChenShonKR,GundersenYang,LiXMGao,YangLZKodai,YangLZ1,YangLZ2}.
 For the proof of Theorem \ref{T1}, we need one of such results as follows.
 \begin{thm}\label{T2}\cite[Theorem 1.1]{LiXMGao}  {\it Let $p(z)$ be a nonconstant polynomial and $Q(z)$ be a nonzero polynomial, then the hyper-order of $f$ is just equal to the degree of $p(z)$.}
\end{thm}
% \begin{thm}\label{T4}\cite[Theorem 2.3]{CaoAus}  {\it Let $p(z)$ be a transcendental entire function of order $\rho(p)\leq 1/2$ and $Q(z)$ be a small entire function respect to $f$, then the hyper-order of $f$ is infinite.}
%\end{thm}
In the above theorem, $Q(z)$ can be a constant. In order to achieve our goal, we shall also study the hyper-order of solutions of a complex differential equation firstly. Following is the statement.

\begin{thm} \label{T3} {\it Suppose $A(z)=-e^{p(z)}$, where $p(z)$ is a transcendental
entire function with nonzero finite order and $B(z)$ is an entire function with nonzero finite hyper-order. If $$\frac{1}{\rho_2(A)}+\frac{1}{\rho_2(B)}>2,$$ then every solution $f(\not\equiv0)$ of equation
\bea\label{1.2}
f^{(k)}+A(z)f'+B(z)f=0,\ \ k(\in\mathbb{N})>1, \eea is of infinite  hyper-order.}\\
\end{thm}
The method of proving this theorem is originally from Rossi \cite{Rossi}. It was also used by Cao \cite{CaoAus} to affirm the Br\"{u}ck conjecture when $f(z)$ is of hyper-order $1/2$.

 \setcounter{section}{1}
\section{Preliminary lemmas}

\begin{lemma}\label{L1}
\cite{GundersenLogarithmic} Let $f$ be a transcendental meromorphic function. Let $\alpha>1$ be a constant, and $k, j$ be integers satisfying $k>j\geq 0$. There exists a set $E\subset [0,2\pi)$ which has zero linear measure, such that if $\theta\in[0,2\pi)\setminus E$, then there is a constant $R(=R(\theta))>0$ such that  \bea\label{2.1} \left|\frac{f^{(k)}(z)}{f^{(j)}(z)}\right|\leq K\left[\frac{T(\alpha r,f)}{r}(\log r)^{\alpha}\log T(\alpha r,f)\right]^{k-j}\eea holds for all $z$ satisfying $\arg z=\theta$ and $|z|\geq R$.
\end{lemma}

 The following Lemma is proved in \cite{Rossi} by using \cite[Theorem \uppercase\expandafter{\romannumeral3}.68]{Tsuji}. Some notations are needed to state it. Suppose $D$ is a domain in $\mathbb{C}$. For each $r\in \mathbb{R}^+$ set $\theta^*_D(r)=\theta^*(r)=+\infty$ if the entire circle $|z|=r$ lies in $D$. Otherwise, let $\theta^*_D(r)=\theta^*(r)$ be the measure of all $\theta$ in $[0,2\pi)$ such that $re^{i\theta}\in D$. As usual, we define the order $\rho(u)$ of a function $u$ subharmonic in the plane as

 $$\rho(u):=\limsup_{R\rightarrow+\infty}\frac{\log M(r,u)}{\log r},$$
here $M(r,u)$ denote the maximum modulus of subharmonic function $u$ on a circle of radius $r$.
\begin{lemma}\label{L2}
\cite{Rossi} Let $u$ be a subharmonic function in $\mathbb{C}$ and let $D$ be an open component of $\{z:u(z)>0\}$. Then
\bea\label{2.2} \rho(u)\geq \limsup_{R\rightarrow+\infty}\frac{\pi}{\log R}\int_1^R\frac{dt}{t\theta^*_D(t)}.\eea
 Furthermore, given $\varepsilon>0$, define $F=\{r:\theta^*_D\leq \varepsilon\pi\}$. Then \bea \label{2.3} \limsup_{R\rightarrow+\infty}\frac{1}{\log R}\int_{F\cap [1,R]}\frac{dt}{t}\leq\varepsilon \rho(u).\eea

\end{lemma}

\begin{lemma}\label{L3}
\cite{Eremenko,Rossi} Let $l_1(t)>0$, $l_2(t)>0(t\geq t_0)$ be two measurable functions on $(0,+\infty)$ with $l_1(t)+l_2(t)\leq (2+\varepsilon)\pi$, where $\varepsilon>0$. If $G\subseteq(0,+\infty)$ is any measurable set and
\bea \pi\int_G\frac{dt}{tl_1(t)}\leq \alpha\int_G\frac{dt}{t},\hskip 0.5cm \alpha\geq \frac{1}{2},\eea
then \bea \pi\int_G\frac{dt}{tl_2(t)}\geq\frac{\alpha}{(2+\varepsilon)\alpha-1}\int_G\frac{dt}{t}.\eea
\end{lemma}

\setcounter{section}{2}
\section{Proof of Theorems}

\noindent{\it Proof of Theorem \ref{T3}}

Suppose that $\rho_2(f)<+\infty$, and would obtain the assertion by reduction to a contradiction.
From Lemma \ref{L1} and the definition of growth order, there exist constants $K>0$, $\beta>1$ and $C=C(\varepsilon)$ (depending on $\varepsilon$) such that
\bea\label{3.1} \left|\frac{f^{(j)}(z)}{f(z)}\right|\leq K\left(\frac{T(\beta r,f)}{\log r}(\log r)^{\beta}\log T(\beta r,f)\right)^j\leq \exp\{r^C\},\hskip 0.5cm j=1,k\eea
holds for all $r>r_0=R(\theta)$ and $\theta\not\in J(r)$, where $J(r)$ is a zero linear measure set. For any given positive small $\varepsilon$, we give $m(J(r))\leq \varepsilon\pi$ .

Fix $\varepsilon>0$ and take a positive integer $N$ such that $ N>C=C(\varepsilon)$.
Since $B(z)$ is an entire function with infinite order and $\log|\log B(z)|=\log\log|B(z)|+o(1)$ as $|B(z)|$ sufficiently large, it's feasible to define the set $\widetilde{D}:=\{z:\log|\log B(z)|-N\log|z|>0\}\cap \{z:\log\log |B(z)|-N\log|z|>0\}$ (here, and in the sequel, taking the principal value of complex logarithm). From \cite{Hayman}, we know that if $u(z)$ is analytic in a domain $D$, then $\log|u(z)|$ is subharmonic in $D$. Since $\frac{\log B(z)}{z^N}$ is analytic in $\widetilde{D}$, the function $\log|\frac{\log B(z)}{z^N}|=\log|\log B(z)|-N\log|z|$ is subharmonic in the open set $\widetilde{D}$. Choose
one unbounded component of $\widetilde{D}$, called $D_1$, such that if we define
\bna u(z)=
\begin{cases}
\log|\log B(z)|-N\log|z|, &z\in D_1, \cr 0, &z\in \mathbb{C}\setminus D_1, \end{cases}
\ena
then $u(z)$ is subharmonic in $\mathbb{C}$ with \bea\label{3.2} \rho(u)\leq\rho_2(B).\eea
Let $D_2$ be an unbounded component of the set $\{z:\log |\log (-e^{-p(z)})|>0\}\cap\{z:\log \log| e^{-p(z)}|>0\}$, such that if we define
\bna v(z)=
\begin{cases}
\log|\log(-e^{-p(z)})|, &z\in D_2, \cr 0, &z\in \mathbb{C}\setminus D_2, \end{cases}
\ena
then $v(z)$ is subharmonic in $\mathbb{C}$ with $\rho(v)=\rho_2(A)$. Moreover, define $D_3:=\{re^{i\theta}:\theta\in J(r)\}$.
For the above given $\varepsilon$, if $(D_1\cap D_2)\setminus D_3$ contains an unbounded sequence $\{r_ne^{i\theta_n}\}$, by \eqref{1.2} and \eqref{3.1} and together with the properties the sets $D_1, D_2$ have we get
\bna \exp\{r_n^N\}<|B(r_ne^{i\theta_n})|\leq \left|\frac{f^{(k)}(r_ne^{i\theta_n}))}{f(r_ne^{i\theta_n})}\right|+\left|{e^{p(r_ne^{i\theta_n})}}\right|\left|\frac{f'(r_ne^{i\theta_n}))}{f(r_ne^{i\theta_n})}\right|
\leq 2\exp\{r_n^{C}\},\ena
this clearly contradicts $N>C$ for $n$ large enough. Thus we could assume that $(D_1\cap D_2)\setminus D_3$ is bounded for arbitrary $\varepsilon$, this implies that for $r\geq r_1\geq r_0$, ($r_0$ is from the bottom of \eqref{3.1}) \bna K_r:=\{\theta: re^{i\theta}\in D_1\cap D_2\}\subseteq J(r).\ena Obviously, \bea\label{3.3} m(K_r)\leq \varepsilon\pi.\eea
(We remark here that if $D_1$ and $D_2$ were disjoint, the proof of Theorem \ref{T3} would follow easily from \eqref{2.2} and Lemma \ref{L3}. In fact, from \eqref{2.2} and \eqref{3.3} we can deduce that the sets  are disjoint in some sense.)
Define

\bna l_j(t)=
\begin{cases}
2\pi, & \rm {if}\ \theta^*_{D_j}(t)=\infty, \cr \theta^*_{D_j}(t), & \rm{otherwise}, \end{cases}
\ena
for $j=1,2$. Since $D_1$ and $D_2$ are unbounded open sets we have that $l_1(t)>0, l_2(t)>0$ for $t$ sufficiently large, and \eqref{3.3} gives \bea\label{3.4} l_1(t)+l_2(t)\leq 2\pi+\varepsilon\pi.\eea Set
\bea\label{3.20} \alpha:=\limsup_{R\rightarrow\infty}\frac{\pi}{\log R}\int_1^R \frac{dt}{tl_1(t)}.\eea
From \eqref{3.20} and the fact $l_1(t)\leq2\pi$, it's clear that \bea \alpha\geq\frac{1}{2\log R}\int_1^R \frac{dt}{t}=\frac{1}{2}.\nonumber\eea
Also by \eqref{3.20} we have \bea \pi\int_1^R \frac{dt}{tl_1(t)}\leq\alpha\log R=\alpha\int_1^R\frac{dt}{t}.\nonumber\eea
Then the conditions of Lemma \ref{L3} are satisfied, we obtain \bea \pi\int_1^R \frac{dt}{tl_2(t)}\geq \frac{\alpha}{(2+\varepsilon)\alpha-1}\int_1^R\frac{dt}{t}=\frac{\alpha}{(2+\varepsilon)\alpha-1}\log R,\nonumber\eea
this means, \bea\label{3.5} \limsup_{R\rightarrow+\infty}\frac{\pi}{\log R}\int_1^R \frac{dt}{tl_2(t)}\geq \frac{\alpha}{(2+\varepsilon)\alpha-1}.\eea
Define the sets \bna B_j:=\{r:\theta^*_{D_j}(r)=+\infty\}\ena
for $j=1,2$. If $r\in B_1$ and $r\geq r_1$, then $\theta^*_{D_2}(r)\leq \varepsilon\pi$ by \eqref{3.4}. Thus $B_1\subseteq \{r: \theta^*_{D_2}(r)\leq \varepsilon\pi\}$. By Lemma \ref{L2} we have
\bea\label{3.6} \limsup_{R\rightarrow\infty}\frac{1}{\log R}\int_{B_1\cap[1,R]}\frac{dt}{t}\leq \varepsilon\rho_2(-e^{-p(z)})=\varepsilon\rho_2(-e^{p(z)})=\varepsilon\rho_2(A).\eea
The equality follows by the first Nevanlinna theorem.
Set $\widetilde{B_j}=R^+\setminus B_j, j=1,2$. Then \eqref{2.2},
\eqref{3.20} and \eqref{3.6} give \bna\rho(u)&\geq&\limsup_{R\rightarrow\infty}\frac{\pi}{\log R}\int_1^R\frac{dt}{t\theta^*_{D_1}(t)}\nonumber\\
&=& \limsup_{R\rightarrow\infty}\frac{\pi}{\log R}\int_{\widetilde{B_1}\cap[1,R]}\frac{dt}{t\theta^*_{D_1}(t)}\nonumber\\
&=& \limsup_{R\rightarrow\infty}\frac{1}{\log R}\left[\pi\int_1^R\frac{dt}{t l_1(t)}-\frac{1}{2}\int_{{B_1}\cap[1,R]}\frac{dt}{t}\right]\nonumber\\
&\geq& \alpha-\frac{\varepsilon\rho_2(A)}{2},
\ena
which together with \eqref{3.2} show \bea\label{3.7}  \rho_2(B)\geq \alpha-\frac{\varepsilon\rho_2(A)}{2}.\eea
For the set $B_2$, we have the similar result as follows by the above arguments for $B_1$. If $r\in B_2$ and $r\geq r_1$, then $\theta^*_{D_1(r)}\leq\varepsilon\pi$. Thus $B_2\subseteq\{r:\theta^*_{D_1(r)}\leq\varepsilon\pi\}$.
 Also from Lemma \ref{L2} we get \bea\label{3.8} \limsup_{R\rightarrow\infty}\frac{1}{\log R}\int_{B_2\cap[1,R]}\frac{dt}{t}\leq \varepsilon\rho(u).\eea
Combining \eqref{2.2}, \eqref{3.5} with \eqref{3.8} we obtain
\bea \label{3.9}\rho_2(A)=\rho(v)&\geq&\limsup_{R\rightarrow\infty}\frac{\pi}{\log R}\int_1^R\frac{dt}{t\theta^*_{D_2}(t)}\nonumber\\
&=& \limsup_{R\rightarrow\infty}\frac{\pi}{\log R}\int_{\widetilde{B_2}\cap[1,R]}\frac{dt}{t\theta^*_{D_2}(t)}\nonumber\\
&=& \limsup_{R\rightarrow\infty}\frac{1}{\log R}\left[\pi\int_1^R\frac{dt}{t l_2(t)}-\frac{1}{2}\int_{{B_2}\cap[1,R]}\frac{dt}{t}\right]\nonumber\\
&\geq& \frac{\alpha}{(2+\varepsilon)\alpha-1}-\frac{\varepsilon\rho(u)}{2}.
\eea
Since $\frac{\alpha}{(2+\varepsilon)\alpha-1}$ is a monotone decreasing function of $\alpha$,
inequalities \eqref{3.2}, \eqref{3.7} and \eqref{3.9} give
\bna\rho_2(A)\geq \frac{ \rho_2(B)+\frac{\varepsilon\rho_2(A)}{2}}{(2+\varepsilon)(\rho_2(B)+\frac{\varepsilon\rho_2(A)}{2})-1}-\frac{\varepsilon\rho_2(B)}{2}.\ena
Note that $\varepsilon$ is arbitrary positive small and $\rho_2(A), \rho_2(B)$ are finite, we obtain \bna \rho_2(A)\geq \frac{\rho_2(B)}{2\rho_2(B)-1}.\ena
It can be transformed into $$ \frac{1}{\rho_2(A)}+\frac{1}{\rho_2(B)}\leq 2,$$
which contradicts the assumption. Thus, every solution $f(\not\equiv0)$ of equation \eqref{1.2} is of infinite hyper-order.

\hskip 0.5cm

 \noindent{\it Proof of Theorem \ref{T1}}

By the assumption that $f$ is a nonconstant entire function with hyper-order $\rho_2(f)<1$, obviously we have \bna \rho_2\left(\frac{f^{(k)}-a}{f-a}\right)<1.\ena
Noting that $f^{(k)}-a$ and $f-a$ share $0$ CM and by the result of the essential
part of the factorization theorem for meromorphic function of finite iterated order \cite[Satz12.4]{JankVolkmann}, we have
\bea\label{3.10} \frac{f^{(k)}-a}{f-a}=e^{p(z)},\eea
 where $p(z)$ is an entire function with $\rho(p(z))=\rho_2(e^{p(z)})<1$.
Suppose that $p(z)$ is not a constant. Set $F:=f-a$, clearly it's not identically equal to zero, then $f^{(k)}=F^{(k)}$. Equation \eqref{3.10} becomes \bea\label{3.11} F^{(k)}-e^{p(z)}F=a.\eea
 Differential both sides we obtain \bea\label{3.12} F^{(k+1)}-e^{p(z)}F'-p'(z)e^{p(z)}F=0.\eea
 Set $A(z)=-e^{p(z)}, B(z)=-p'(z)e^{p(z)}$, thus $\rho(p)=\rho_2(A)=\rho_2(B)<1$. If $p(z)$ is a nonconstant polynomial, applying Theorem \ref{T2} to equation \eqref{3.11} we deduce that $\rho_2(F)=\rho_2(f)$ is equal to a positive integer, which contradicts the assumption $\rho_2(f)<1$.
 If $p(z)$ is transcendental with $\rho(p)<1/2$, by the proof of Theorem 2 Case (3) in \cite{ChenShonKR}, it follows a contradiction. If $p(z)$ is transcendental with $1/2\leq\rho(p)<1$, applying Theorem \ref{T3} to equation \eqref{3.12}, we have $\rho_2(F)=\rho_2(f)$ is infinite, also contradicts the assumption $\rho_2(f)<1$. Therefore, $p(z)$ must be a constant, and thus $e^{p(z)}$ is just a nonzero constant. Then, we complete the proof.

 \vskip 1cm
%\vskip 3cm
%\noindent{\bf Competing interests.}  The author declares that there
%is no conflict of interests regarding the publication of this
%article.\\
%%
%\noindent{\bf Author's contributions.} The author completed the
%paper and approved the final manuscript.\\
\noindent{\bf Acknowledgements.}
 This work was supported by NSFC (no. 11426035). The author wishes to thank Prof. Tingbin Cao of Nanchang University for his kind help and Post Doctor Xiao Yao of Fudan University for helpful discussion to make it more readable.

% and the Key Project
%of Natural Science Foundation of Educational Committee of Henan
%Province under Grant no. 15A110008.
%The authors wish to express their thanks to the referee for his/her
%valuable suggestions and comments.


\begin{thebibliography}{16}

\bibitem{Bruck}{}R. Br\"{u}ck, On entire functions which share one value CM with their first derivative, Results in
Math. 30(1996), 21-24.

%
\bibitem{CaoPolon} T. B. Cao, Growth of solutions of a class of complex differential equations, Ann. Polon. Math.
95(2009), No. 2, 141-152.
\bibitem{CaoAus} T. B. Cao, On the Br\"{u}ck conjecture, Bull. Aust. Math. Soc. 93 (2016), no. 2, 248--259.
%\bibitem{Chen} Z.-X. Chen, K.-H. Shon, The growth of solutions of differential equations with
%coefficients of small growth in the disc, J. Math. Anal. Appl.
%297(2004) 285-304.



\bibitem{Changjm}{}
J. M. Chang and Y. Z. Zhu, Entire functions that share a small function with their derivatives, J.
Math. Anal. Appl. 351(2009), 491-496.
\bibitem{ChenShonTW}{}
Z. X. Chen and K. H. Shon, On conjecture of R. Br\"{u}ck concerning the entire function sharing one
value CM with its derivative, Taiwanese J. Math. 8(2004), No. 2, 235-244.
\bibitem{ChenShonKR}{}
Z. X. Chen and K. H. Shon, On the entire function sharing one value CM with k-th derivatives, J.
Korean Math. Soc. 42(2005), No. 1, 85-99.
%

\bibitem{ChenZhang}{}
Z. X. Chen and Z. L. Zhang, Entire functions sharing fixed points with their higher order
derivatives, Acta Math Sin. Chinese Ser. 50(2007), No. 6, 1213-1222. (in Chinese)

\bibitem{Eremenko}{}
A. Eremenko, Growth of entire and meromorphic functions on asymptotic curves, Sibirsk Mat.Zh. 21 (1980),39-51; English transl. in Siberian Math. J. (1981), 673-683.
\bibitem{GundersenLogarithmic} G. Gundersen, Estimates for the logarithmic derivative of a meromorphic function, Plus Similar
Estimates, J. London Math. Soc. 37(1988), No.2, 88-104.
\bibitem{GundersenYang} G. G. Gundersen and L. Z. Yang, Entire functions that share one value with one or two of their derivatives, J. Math. Anal. Appl. 223(1998), 88-95.
\bibitem{Hayman}W. Hayman, P. Kennedy, Subharmonic functions, Academic press, London, 1976.
\bibitem{JankVolkmann}G. Jank and L. Volkmann, Meromorphe Funktionen und Differentialgleichungen, Birk\"{a}user, 1985.
%\bibitem{Gundersen}G. Gundersen, Finite order solutions of second order linear differential equations, Trans.
%Amer. Math. Soc. 305 (1988), 415-429.

\bibitem{Laine}I. Laine, Nevanlinna theory and complex differential equations, Walter de Gruyter,
Berlin, 1993.
\bibitem{LiXiaomin} X. M. Li, An entire function and its derivatives sharing a polynomial, J. Math. Anal. Appl.
330(2007), 66-79.
%\bibitem{HellersteinRossi}{}
%S. Hellerstein, J. Miles, J. Rossi, On the growth of solutions $f''+gf'+hf=0$, Trans. Amer. Math. Soc. 324(1991), 693-706.

%\bibitem{Hille}{}
%E. Hille, Ordinary differential equations in the complex domain, New York, Wiley, 1976.
%\bibitem{Holland}A. S. B. Holland, Introduction to the theory of entire functions, Academic Press, New York-London, 1973.
\bibitem{LiXMGao}X. M. Li and C. C. Gao, Entire functions sharing one polynomial with their derivatives, Proc.
Indian Acad. Sci. (Math. Sci.) 118(2008), No. 1, 13-26.
%\bibitem{Ishizaki}K. Ishizaki and K. Tohge, On the complex oscillation of some linear
%differential equations, J. Math. Anal. Appl. 206(1997), 503-517.

\bibitem{RubelYang}L. A. Rubel and C. C. Yang, Values shared by an entire function and its derivative, Complex
Analysis, Lecture Notes in Mathematics 599(1977), 101-103


%
%

\bibitem{Rossi}J. Rossi, Second order differential equations with transcendental coefficients, Proc. Amer. Math. Soc. 97(1986), 61-66.

\bibitem{Tsuji}M. Tsuji, Potential theory in modern function theory, Maruzen, Tokyo, 1959.
\bibitem{WangJun}J. Wang, Uniqueness of entire function sharing a small function with its derivative, J. Math. Anal.
Appl. 362(2010), No. 2, 387-392.
\bibitem{YangLZKodai}L. Z. Yang, Solutions of a differential equation and its applications, Kodai Math. J. 22(1999), No.
3, 458-464

\bibitem{YangLZ1}L. Z. Yang, Entire functions that share one value with one of their derivatives, in Finite or Infinite
Dimensional Complex Analysis, Fukuoka 1999 (Proceedings of a Conference), Lecture Notes in
Pure and Applied Mathematics, Vol. 214, Marcel Dekker, New York, 2000, pp. 617-624.
\bibitem{YangLZ2}L. Z. Yang, The growth of linear differential equations and their applications, Israel J. Math.
147(2005), 359-370.

%\bibitem{Wushengjian}S.-J. Wu, On the growth of solutions of second order linear differential equations
%in an angle, Complex Variables and Elliptic Equations, 1994(24),
%241-248.
%
%\bibitem{XuJunfeng}J.-F. Xu, H.-X. Yi, Solutions of higher order linear differential equations in an angle,
%Applied Mathematics Letters, 2009, 22(4), 484-489.
%\bibitem{WuPCZhuJ}P. C. Wu, J. Zhu, On the growth of solutions to
%the complex differential equation $f''+Af'+Bf=0$, Sci. China
%Math., 2011, 54(5): 939-947.
\bibitem{YY} C. C. Yang, H. X. Yi, Uniqueness theory of meromorphic functions, Kluwer Academic Publishers, Dordrecht,
2003.

%\bibitem{ZhangGH}
%G.-H. Zhang, Theory of entire and meromorphic functions, deficient
%and asymptotic values and singular directions, Translations of
%Mathematical Monographs, {Amer. Math. Soc. }122(1993).
%\bibitem{ZhengJH2}
%J.-H. Zheng, Value distribution of meromorphic functions, Tsinghua
%University Press, Beijing and Springer-Verlag Berlin Heidelberg,
%2010.
%\bibitem{ZhengJH3}
%J. H. Zheng, S. Wang, Z. Huang, Some properties of Fatou and Julia
%sets of transcendental meromorphic functions, \emph{Bull. Aust.
%Math. Soc.} 66 (2002) 1-8.
\end{thebibliography}
\end{document}